\documentclass[reqno,oneside,10pt]{amsart}

\usepackage{graphicx}
\usepackage{enumerate}
\usepackage{xy}
\usepackage{amsmath,amsfonts,amssymb}
\usepackage[latin9]{inputenc}
\newcommand{\N}{{\mathbb N}}

\newcommand{\prho}{\overline{\rho}} 

\newtheorem{lemma}{Lemma} 
\newtheorem{prop}{Proposition} 
\newtheorem{cor}{Corollary} 
\newtheorem{teo}{Theorem}

\newtheorem{defi}{Definition}

\newcommand{\vai}{\rightarrow}

\newcommand{\ee}{{\mathbf e}}

\newcommand{\hh}{{\mathbf h}}

\newcommand{\um }{1_A}
\usepackage{latexsym}

\usepackage{xspace}
\usepackage{amscd}
\usepackage{amssymb}
\usepackage{amsfonts}

\newcommand{\lacts}{\rhd}
\newcommand{\hits}{\rhd}

\title[Enveloping Actions for Partial Hopf Actions]{Enveloping Actions for Partial Hopf Actions}

\author[M.M.S. Alves]{Marcelo \ Muniz \ S. \ Alves}
\address{Departamento de Matem\'atica, Universidade Federal do Paran\'a, Brazil}
\email{marcelo@mat.ufpr.br}
\author[E. Batista]{Eliezer Batista}
\address{Departamento de Matem\'atica, Universidade Federal de Santa Catarina, Brazil}
\email{ebatista@mtm.ufsc.br}

\begin{document}

\begin{abstract}
 Motivated by partial group actions on unital algebras, in this
article we extend many results obtained by Exel and Dokuchaev \cite{dok}  to the context of partial
actions of Hopf algebras, according to Caenepeel and Jansen \cite{caen06}. First, we generalize the
theorem about the existence of an enveloping action, also known as the
globalization theorem. Second, we construct a Morita context between the
partial smash product defined by the authors of \cite{caen06} and the smash product related
to the enveloping action. Third, we dualize the globalization theorem to
partial coactions and finally, we define partial representations of Hopf
algebras and show some results relating partial actions and partial
representations.
\end{abstract} 

\maketitle

\section{Introduction}
Partial group actions were first defined by R. Exel in the context of operator algebras and
they turned out to be a powerful tool in the study of $C^*$-algebras generated
by partial isometries on a Hilbert space \cite{ruy}. The developments
originated by the definition of partial group actions include crossed products
\cite{quigg}, partial representations \cite{ruy2,dok1} and soon this theme
became an independent topic of interest in ring theory
\cite{dok,ferrero}. Now, the results are formulated in a purely 
algebraic way, independent of the  $C^*$ algebraic techniques which originated 
them. 

A partial action $\alpha$ of a group $G$ on a (possibly non-unital) $k$-algebra $A$ is a pair of families of sets and maps indexed by $G$, $\alpha = (\{\alpha_g\}_{g \in G}, \{D_g\}_{g \in G})$, where each $D_g$ is an ideal of $A$ and each $\alpha_g$ is an algebra isomorphism $\alpha: D_{g^{-1}} \rightarrow D_g$ satisfying the following conditions:
\begin{enumerate}
\item[(i)] $D_e = A$ and $\alpha_e = I_A$;
\item[(ii)] $\alpha_g(D_{g^{-1}} \cap D_h) = D_g \cap D_{gh}$ for every $g,h \in G$;
\item[(iii)] $\alpha_g(\alpha_h(x)) = \alpha_{gh}(x)$ for every 
$x \in D_{g^{-1}} \cap D_{(gh)^{-1}}.$ 
\end{enumerate}

A first example of partial action is the following: If $G$ acts on a 
algebra $B$ by automorphisms and $A$ is an ideal of $B$, then we
have a partial action $\alpha$ on $A$ in the following manner: letting
$\beta_g$ stand for the automorphism corresponding to $g$, take 
$D_g = A \cap \beta_g(A)$, and define $\alpha_g: D_{g^{-1}} \vai D_g$ 
as the restriction of the automorphism $\beta_g $ to $D_g$. 

Partial Hopf actions were motivated by an attempt to
generalize the notion of partial Galois extensions of commutative rings, first
introduced by M. Dokuchaev, M. Ferrero and A. Paques in \cite{paques}. The
first ideas towards partial Hopf actions were introduced by S. Caenepeel 
and E. de Groot in \cite{caen04}, using the concept of Galois
coring. Afterwards, S. Caenepeel and K. Janssen defined partial actions and
partial co-actions of a Hopf algebra $H$ on a unital algebra $A$ using the
notions of partial entwining structures \cite{caen06}; in particular, partial actions of $G$ determine 
partial actions of the group algebra $kG$ in a natural way. In the
same article, the authors also introduced the concept of partial 
smash product, which in the case of the group algebra $kG$, turns 
out to be the crossed product by a partial action 
$A\rtimes_{\alpha} G$. Further developments in the theory of
partial Hopf actions were done by C. Lomp in \cite{lomp07}, where the author
pushed forward classical results of Hopf algebras concerning smash products, 
like the Blattner-Montgomery and Cohen-Montgomery theorems \cite{susan}.

Certainly, the theory of partial actions of Hopf algebras remains as a huge
landscape to be explored, and this present work intends to generalize some
results for partial group actions into the context of partial Hopf actions. We
divided this paper as follows:

In section 2, we prove the theorem of existence of an enveloping action for a partial Hopf action, i.e, we prove that  if $H$ is a Hopf algebra which acts partially on a unital algebra $A$, 
then there exists an $H$ module algebra $B$ such that $A$ is isomorphic 
to a right ideal of $B$, and the restriction of the action of $H$ to this ideal is 
equivalent to the partial action of $H$ on $A$. The {\it uniqueness} of the enveloping action is treated separately;  we introduce the concept of minimal enveloping action and prove that the existence and uniqueness of such an action for every partial action. The question for enveloping actions for partial group actions arises naturally when we
consider the basic example of partial action as a restriction of a global
action of a group $G$ on an algebra $B$ to an ideal $A\trianglelefteq B$. What
conditions on the partial action enables us to say that this partial action
is a restriction of a global action? The first result concerning enveloping
actions was proved in the context of $C^*$ algebras in \cite{abadie}; to this
intent, the author used techniques of Fell Bundles and Hilbert $C^*$ modules. 
A purely algebraic version of this theorem on enveloping actions only 
appeared in \cite{dok}. Basically, the same ideas for the proof in the 
group case are present in the Hopf algebraic case as we shall see later.

In section 3, we show the existence of a Morita context between the partial
smash product $\underline{A\# H}$, where $H$ is a Hopf algebra which acts
partially on the unital algebra $A$, and the smash product $B\# H$, where $B$
is an enveloping action of $A$. This result can also be found in \cite{dok}
for the context of partial group actions.

In section 4, we discuss the existence of an enveloping co-action of a Hopf
algebra $H$ on a unital algebra $A$. There, we dualize this partial co-action
of $H$ to a partial action of $H^*$ (in fact, the finite dual $H^\circ$), we take an enveloping action and then check whether the $H^\circ$ module $B$ of the
enveloping action is a rational module. If this occurs, one
dualizes again to obtain
a structure of $H$ comodule algebra in $B$; this is our enveloping co-action.

In section 5, we introduce the notion of partial representation of a Hopf
algebra. We show that, under certain conditions on the algebra $H$, the
partial smash product $\underline{A\# H}$ carries a partial representation of
$H$. 

\section{Enveloping actions}

\subsection{Partial Hopf actions}

We recall that a left action of a Hopf algebra $H$ on an algebra $A$ is  
a linear mapping $\alpha: H \otimes A \vai A$, which we will denote by 
$\alpha(h \otimes a) = h \hits a$, such that  
\newpage

\begin{enumerate}
\item[(i)] $h \hits (ab) = \sum (h_{(1)} \hits a) (h_{(2)} \hits b),$
\item[(ii)] $1 \hits a =a$
\item[(iii)]  $h \hits (k \hits a) = hk \hits a$ 
\item[(iv)] $h \hits \um = \epsilon(h) \um. $
\end{enumerate}
We also say that $A$ is an $H$ {\it module algebra}. Note that (ii) and (iii) say that $A$ is a left $H$-module.

In \cite{caen06}, Caenepeel and Jansen defined a weaker version of an action, 
called a partial action. A partial action of the Hopf algebra $H$ on the 
algebra $A$ is a linear mapping $\alpha: H \otimes A \vai A$, denoted 
here by $\alpha(h \otimes a) = h \cdot a$, such that  
\begin{enumerate}
\item[(i)] $h \cdot (ab) = \sum (h_{(1)} \cdot a) (h_{(2)} \cdot b),$
\item[(ii)] $1 \cdot a = a,$
\item[(iii)] $h \cdot (g \cdot a) = \sum (h_{(1)} \cdot \um) ((h_{(2)} g) \cdot a).$
\end{enumerate}

We will also call $A$ a {\it partial} $H$ {\it module algebra}. 
It is easy to see that every action is also a partial action.

As a basic example, consider  a partial action $\alpha$ of a group $G$ on an
unital algebra $A$. Suppose that each $D_g$ is also a unital algebra, that is,
$D_g$ is of the form $D_g =A1_g$ then there is a partial action of the group
algebra $kG$ on $A$ defined on the elements of the basis by
\begin{equation}
\label{actiongroupalgebra}
g\cdot a =\alpha_g (a 1_{g^{-1}}),
\end{equation}
and extended linearly to all elements of $kG$. In order to see that this
action satisfies the relations (i), (ii) and (iii) of the definition of partial action above, let us remember some
facts about the partial action $\alpha$. First, the elements $1_g \in D_g$ are
central idempotents in the algebra $A$ and are given by $1_g = g \cdot \um$, second, the unity of the ideal 
$D_g \cap D_h$ is the product $1_g 1_h$ and finally, since 
$\alpha_g (D_{g^{-1}} \cap D_h )= D_g \cap D_{gh}$ and each $\alpha_g$ is an
isomorphism, we have $\alpha_g (1_{g^{-1}} 1_h)=1_g 1_{gh}$. Then, the action
(\ref{actiongroupalgebra}) satisfies:
\begin{eqnarray}
g\cdot (ab) &=& \alpha_g (ab1_{g^{-1}})= 
\alpha_g (a1_{g^{-1}} b1_{g^{-1}})=\nonumber\\
&=& \alpha_g (a1_{g^{-1}})\alpha_g (b1_{g^{-1}})=
(g\cdot a)(g\cdot b),\nonumber
\end{eqnarray}
\begin{eqnarray} 
e\cdot a &=& \alpha_e (a1_{e^{-1}})= I_A (a \um)= a, \nonumber
\end{eqnarray}
\begin{eqnarray}
h\cdot (g\cdot a) &=& \alpha_h (\alpha_g (a1_{g^{-1}})1_{h^{-1}})=\nonumber\\
&=& \alpha_h (\alpha_g (a1_{g^{-1}})1_g 1_{h^{-1}})=\nonumber\\
&=& \alpha_h (\alpha_g (a1_{g^{-1}})\alpha_g( 1_{g^{-1}} 1_{g^{-1} h^{-1}}))=\nonumber\\
&=& \alpha_h (\alpha_g (a1_{g^{-1}} 1_{g^{-1} h^{-1}}))=\nonumber\\
&=& \alpha_{hg} (a1_{g^{-1}} 1_{g^{-1} h^{-1}})=\nonumber\\
&=& \alpha_{hg} (a1_{g^{-1} h^{-1}})\alpha_{hg}(1_{g^{-1}} 1_{g^{-1} h^{-1}})=\nonumber\\
&=& \alpha_{hg} (a1_{g^{-1} h^{-1}})1_{hg} 1_h=1_h \alpha_{hg} (a1_{g^{-1} h^{-1}})\nonumber\\
&=& \alpha_h (\um 1_{h^{-1}}) \alpha_{hg} (a1_{g^{-1} h^{-1}})=
(h\cdot \um)(hg\cdot a).\nonumber
\end{eqnarray}

Note that we have also proved that $h \cdot (g \cdot a) = (hg\cdot a) 1_h = (hg\cdot a)(h\cdot \um)$. In general, it is not true that any partial action of $kG$ induces automatically a partial group action of $G$. We mention that in \cite{caen06} the authors consider a slight generalization of partial group actions, where the idempotents $1_g$ are not necessarily central and  $D_g$ is the right ideal $D_g = 1_g A$; in this case, it can be proven that there is a bijective correspondence between partial group actions and partial $kG$-actions on $A$.
 
\subsection{Induced partial actions \label{right}}

There is an important class of examples of
partial Hopf actions induced by total actions. This idea is motivated
by the construction of a partial group action induced by a global action of a
group $G$ on an algebra $B$ by automorphisms.  

Let  $\beta: G \times B \vai B$ be an action of the group $G$ on the 
algebra $B$ by automorphisms, and let $A$ be an ideal of $B$ generated by a central 
idempotent $\um $. Define $D_g = A \cap \beta_g(A)$; then $D_g$ is the ideal 
generated by the central idempotent $1_g = \um  \beta_g(\um )$.  

The partial action $\alpha = (\{\alpha_g\}, \{D_g\})$ induced by $\beta$ on
$A$ is 
\[
\alpha_g(a) = \beta_g(a) \mbox{ for $g \in G$ and $a \in D_{g^{-1}}.$}
\]
This corresponds to a partial action of $kG$ on $A$, given by 
\[
g \cdot a = \alpha_g(a 1_{g^{-1}}).
\]  
Since 
\[
\alpha_g(a 1_{g^{-1}}) = \beta_g( 1_{g^{-1}} a) = 
\beta_g(\um \beta_{g^{-1}}(\um )) \beta_g(a ) = 
\beta_g(\um ) \um \beta_g (a)   = \um \beta_g(a) 
\]
one could also define the partial action by  
$g \cdot a = \um  \beta_g (a)$ (or $g \cdot a = \beta_g (a) \um $). This provides the idea for constructing induced 
partial actions in the Hopf case. 

\begin{prop} \label{induced}
Let $H$ be a Hopf algebra which acts on the algebra $B$, and let $A$ be a
right ideal of $B$ with unity $\um $. 
Then $H$ acts partially on $A$ by 
\[
h \cdot a = \um  (h \hits a) 
\]
\end{prop}
\begin{proof}
The first property is immediate. For the  third, given $h,k \in H$ and $a \in A$, 
\[
\begin{array}{rcl}
h \cdot (k \cdot a) & =& \um (h \hits (\um (k \hits a ))) \\
 & = & \um [\sum (h_{(1)} \hits \um ) (h_{(2)} \hits (k \hits a ))]\\
& = & \sum \um (h_{(1)} \hits \um ) ((h_{(2)} k) \hits a )) = (\ast) \\
\end{array}
\]
and since $\um  (h_{(1)} \hits \um )  \in A $, it follows that $\um (h_{(1)} \hits \um )   = \um  (h_{(1)} \hits \um )  \um$; therefore 
\[
(\ast) =   \sum \um  (h_{(1)} \hits \um ) \um  ((h_{(2)} k) \hits a )) =  \sum (h_{(1)} \cdot \um )  ((h_{(2)} k) \cdot a )).
\]
The second property is proved in an analogous manner.
\end{proof}

We say that the partial action $h \cdot a = \um  (h \hits a) 
$ is the {\it partial action induced by } $B$. 
We remark that  in \cite{caen06} the authors introduce a slightly more general concept of partial group action where the domains $D_g$ are already taken as right ideals.

Although this proposition provides a method for constructing examples, it comes as a surprise that,  in some cases, the induced partial action is total. As we have seen, every partial group action induces a partial $kG$ action, and it is easy to define partial group actions that are not total actions; on the other hand, we prove below that there are no properly partial Hopf actions by universal enveloping algebras of Lie algebras.

\begin{prop} Every induced partial action of an universal 
enveloping algebra $\mathcal{U}(\mathfrak{g})$, with $\mathfrak{g}$ being a
Lie algebra, is in fact a total action.
\end{prop}

\begin{proof} Let $\mathcal{U}(\mathfrak{g})$ be the universal enveloping
  algebra of the Lie algebra $\mathfrak{g}$ acting on an algebra $B$. This
  means in particular that every element of $\mathfrak{g}$ acts as a derivation  in $B$. Let $A$ be a right ideal with unity $\um $, and let 
$X\in \mathfrak{g}$. Then, we have
\[
X\triangleright \um  = X\triangleright \um ^2 =\um  
(X\triangleright \um ) + (X\triangleright \um )\um  .
\]
Using the definition of induced partial action, we conclude that 
the partial action of $X$ on $\um $ is 
\[
X\cdot \um  =\um (X\triangleright \um )= 
\um (X\triangleright \um )\um +
\um (X\triangleright \um ) , 
\]
and since $\um b = \um b \um$ for all $b \in B$, this leads to 
\[
(X\cdot \um )\um =0 \qquad \Rightarrow \quad X\cdot \um  =0.
\]
By a simple induction argument, one can conclude that for every element \linebreak
$\xi = X_1 X_2 \ldots X_n \in \mathcal{U}(\mathfrak{g})$, with 
$X_i \in \mathfrak{g}$, we have $\xi \cdot \um  =0$.

Now, let $\xi$ and $\eta$ two elements of 
$\mathcal{U}(\mathfrak{g})$, with $\xi$ being a monomial of the form 
$X_1 X_2 \ldots X_n$, and let $a\in A$, then we have the partial action
\begin{eqnarray}
\label{ug}
\xi \cdot (\eta \cdot a) &=& \sum (\xi_{(1)} \cdot \um )
((\xi_{(2)}\eta)\cdot a)=\nonumber\\
&=& (1\cdot \um )((\xi \eta) \cdot a) + \nonumber\\
&+& \sum_{k=1}^n \sum_{s\in S_{n,k}} ((X_{s(1)}\ldots X_{s(k)})\cdot \um )
((X_{s(k+1)}\ldots X_{s(n)}\eta)\cdot a), 
\end{eqnarray}
where the sum goes over the $(n,k)$ shuffles for every $1\leq k\leq n$ (we recall that a
$(n,k)$ shuffle is a permutation $s\in S_n$ such that 
$s(1)< \cdots < s(k)$ and \linebreak $s(k+1) < \cdots < s(n)$, and that $S_{n,k}$ denotes the subgroup of $(n,k)$ shuffles). The only nonvanishing term in (\ref{ug}) is
the first, because the terms involving shuffles have monomials of degree greater or
equal to one acting on the unit $\um $ of $A$, which we have already proved
that vanish. Therefore
\[
\xi \cdot (\eta \cdot a)=(1\cdot \um )((\xi \eta) \cdot a)=
(\xi \eta) \cdot a,
\]
which proves that this partial action, is, in fact a total action.
\end{proof}

As another example of an induced partial action, let $H_4$ be the Sweedler 
4-dimensional Hopf algebra, with $\beta = \{1,g,x,xg \}$ 
as basis over the field $k$, where $char(k) \neq 2$. The algebra structure 
is determined by the relations 
\[
g^2 = 1, \qquad  x^2 =0 \qquad \mbox{and} \quad xg = -gx .
\]
The coalgebra structure is given by the  coproducts 
\[
\Delta(g) = g \otimes g, \qquad \Delta(x) = x\otimes1 + g\otimes x ,
\]
and counit $\epsilon(g)=1$, $\epsilon(x)=0$. The antipode $S$ in $H_4$ reads 
\[
S(g)=g , \qquad \mbox{and}\quad S(x) = xg.
\] 
A more suitable basis for the study of ideals of $H_4$ consists of the vectors \linebreak $e_1 = (1+g)/2$, $e_2 = (1-g)/2$, $h_1 = xe_1$, 
$h_2 = xe_2$. The multiplication table of $H_4$ in this new basis elements reads 
\[
\begin{array}{|c|c|c|c|c|}
\hline
 & e_1 & e_2 & h_1 & h_2 \\
\hline 
e_1 & e_1 & 0 & 0 & h_2 \\
\hline 
e_2 & 0 & e_2 & h_1 & 0  \\
\hline 
h_1 & h_1 & 0 & 0 & 0 \\
\hline 
h_2 & 0 & h_2 & 0 & 0 \\
\hline 
\end{array}
\]
The expressions for the coproducts of this new basis, are 
\begin{eqnarray}
\Delta (e_1) &=& e_1 \otimes e-1 + e_2 \otimes e_2 , \nonumber\\
\Delta (e_2) &=& e_1 \otimes e_2 + e_2 \otimes e_1 , \nonumber\\
\Delta (h_1) &=& e_1 \otimes h_1 - e_2 \otimes h_2 + h_1 \otimes e_1 
+h_2 \otimes e_2 , \nonumber\\
\Delta (h_2) &=& e_1 \otimes h_2 - e_2 \otimes h_1 + h_1 \otimes e_2 
+h_2 \otimes e_1 .\nonumber
\end{eqnarray}
The counit calculated in the elements of this new basis take the values
$\epsilon(e_1 ) =1$ and $\epsilon (e_2 )=\epsilon (h_1 )=
\epsilon (h_2 )=0$. Finally, the antipode  of this elements are given by
\[
S(e_1 )=e_1, \quad S(e_2 )=e_2, \quad S(h_1 )=-h_2,   \quad  
S(h_2)=h_1 .
\]

The Hopf algebra $H_4$ acts on itself in the canonical way by the left 
adjoint action, i.e, $h \hits k = \sum_{(h)} h_{(1)} k S(h_{(2)})$. Its action 
is summed up in the table below. 
\[
\begin{array}{|c|c|c|c|c|}
\hline
\hits & e_1 & e_2 & h_1 & h_2 \\
\hline 
e_1 & e_1 & e_2 & 0 & 0 \\
\hline 
e_2 & 0 & 0 & h_1 & h_2 \\
\hline 
h_1 & h_1 - h_2 & h_2 - h_1 & 0 & 0 \\
\hline 
h_2 & 0 & 0 & 0 & 0 \\
\hline 
\end{array}
\]

If $X \subset H_4$, denote by $\langle X \rangle$ the $k$-subspace 
generated by $X$. As one sees directly from the multiplication table, $e_1 H_4 = \langle e_1,h_2 \rangle$ and, since $e_1$ and $h_2$ do not commute, we have to kill the latter in order to get a right ideal with unity. But $\langle h_2 \rangle$  is an ideal of $H_4$ which, given 
the nature of the action, is also a $H_4$-submodule. Hence 
$\overline{B} = H_4/\langle h_2 \rangle$ is a $H_4$-module algebra. In what follows, we denote $x + \langle h_2 \rangle \in \overline{B}$ by $\overline{x}$.

The map 
\[
a\overline{e_1} + b\overline{e_2} + c\overline{h_1} \mapsto 
\left[
\begin{array}{cc} 
a & 0 \\
c & b 
\end{array}
\right]
\]
is an algebra isomorphism. Now, the action of $H_4$ on $\overline{B}$ is as follows:

\[
\begin{array}{|c|c|c|c|}
\hline
\hits & \overline{e_1} & \overline{e_2} & \overline{h_1}  \\
\hline 
e_1 & \overline{e_1} & \overline{e_2} & 0 \\
\hline 
e_2 & 0 & 0 & \overline{h_1}  \\
\hline 
h_1 & \overline{h_1} &  - \overline{h_1} & 0 \\
\hline 
h_2 & 0 & 0 & 0 \\
\hline 
\end{array}
\]

The subspace $A = \langle \overline{e_1} \rangle$ is a right ideal with 
unity in $\overline{B}$. Hence, we have a partial action on $A$ induced 
by the action on $\overline{B}$. This partial action is given by 
\[
e_1 \cdot \overline{e_1} = \overline{e_1}, \qquad e_2 \cdot \overline{e_1} = h_1 \cdot \overline{e_1} = h_2 \cdot \overline{e_1} = 0.\]

Once again, it is easy to see that this partial action is in fact total. This
happens because the subspace $J = \langle e_2, h_1, h_2 \rangle $ is an ideal 
of $H_4$, and hence an $H_4$-submodule of $H_4$ by the left adjoint action; 
therefore $H_4/J$ is an $H_4$-module algebra. Since 
$H_4 = \langle e_1 \rangle \oplus J$ as a vector space, the projection 
of $H_4$ onto $H_4/J$ induces an isomorphism of $H$-module algebras $A \equiv H_4/J$. If one looks 
at the action of $H_4$ on $A$, one gets the same table as for the partial 
action of $H_4$ on $H_4/J$ (via the natural identification of $e_1 + J$ 
with $e_1 + \langle h_2 \rangle$).

We shall prove now that every partial action is induced. 

\subsection{Enveloping actions}

In the context of partial group actions, a natural question arises:
under which conditions can a partial action of a group $G$ on an algebra $A$
be obtained, up to equivalence, from a suitable restriction of a group
action of $G$ on an algebra $B$? In other words, given a partial action
$\alpha=\{ \{\alpha_g \}_{g\in G} \{ D_g\}_{g\in G} \}$ of $G$ on $A$, we
want to extend the isomorphisms $\alpha_g :D_{g^{-1}} \rightarrow D_g$ to automorphisms 
$\beta_g: B\rightarrow B$ of an algebra $B$, such that  $A$ is a subalgebra of $B$ (in fact, an
ideal) and such that this extension is the smallest possible (more precisely,
we impose that $B=\cup_{g\in G} \beta_g (A)$). In this case,
the partial action is said to be an {\it admissible restriction}. We say that an
action $\beta$ of $G$ on $B$ is an enveloping action of a partial action
$\alpha$ of $G$ on $A$ if $\alpha$ is equivalent to an admissible restriction
of $\beta$ to an ideal of $B$.

In the context of partial group actions, it is proved that a partial action
$\alpha$ of a group $G$ on a unital algebra $A$ admits an enveloping action if,
and only if, each of the ideals $D_g \trianglelefteq A$ is a unital
algebra. Moreover, if it exists, this enveloping algebra is unique up to
equivalence (see \cite{dok} theorem 4.5). This is the result we generalize 
here in the context of partial actions of Hopf algebras. 

\begin{defi} Let $A$ and $B$ be two partial $H$-module algebras. We 
will say that a morphism of algebras $\theta: A \vai B$ is a morphism of 
partial $H$-module algebras if $\theta(h \cdot a) = h \cdot \theta (a)$ 
for all $h \in H$ and all $a \in A$. If $\theta$ is an isomorphism,
we say that the partial actions are equivalent.
\end{defi}

\begin{defi} Let $B$ be an $H$-module algebra and let $A$ be a right ideal of $B$ with unity $\um$. We will say that the induced partial action on $A$ is admissible  if $B = H \hits A$.
\end{defi}

\begin{defi}  Let  $A$ be a partial $H$-module algebra. An enveloping action for $A$ is a pair $(B,\theta)$, where 
\begin{enumerate}
\item $B$ is a  (not necessarily unital) $H$-module algebra.
\item The map $\theta: A \vai B$ is a monomorphism of algebras.
\item The sub-algebra $\theta(A)$ is a right ideal in $B$.
\item The partial action on $A$ is equivalent to the induced partial action on $\theta (A)$.
\item  The induced partial action on $\theta(A)$ is admissible.
\end{enumerate}
\end{defi}

We will show now that every partial $H$-action  has an enveloping action.
In \cite{dok}, the authors consider the algebra $\mathcal{F}(G,A)$ of functions from $G$ to $A$. Since there is a canonical algebra monomorphism from $\mathcal{F}(G,A)$ into $Hom_k(kG,A)$, it is reasonable to consider, in the Hopf case, the algebra $Hom_k(H,A)$ in place of $\mathcal{F}(G,A)$. We remind the reader that the product in $Hom_k(H,A)$ is
the convolution product $(f\ast g) (h)=\sum f(h_{(1)}) g(h_{(2)})$, and that $H$ acts on this algebra on the left by 
\[
(h \lacts f)(k) = f(kh)
\]  
where $h, k \in H$ and $f \in Hom_k(H,A)$. 

\begin{lemma}\label{tecnico} Let $\varphi: A \vai Hom_k(H,A)$ be the map given by $\varphi(a)(k) = k \cdot a$. \\
(i) $\varphi$ is a linear injective  map and an algebra morphism.\\
(ii) If $h \in H$ and $a \in A$ then $ \varphi (\um) \ast (h\triangleright \varphi(a)) = \varphi (h\cdot a))$\\
(iii)  If $h \in H$ and $a, b \in A$ then $ \varphi (b) \ast (h\triangleright \varphi(a)) = \varphi (b (h\cdot a))$.
\end{lemma}
\begin{proof} 
It is easy to see that $\varphi$ is linear, because the partial action is bilinear; since $\varphi(a) (1_H) = 1_H \cdot a =a$, it follows that it is also injective. Take $a,b\in A$ and $h\in H$, then we have
\begin{eqnarray}
\varphi (ab)(h) &=& h\cdot (ab)= \sum (h_{(1)} \cdot a)(h_{(2)} \cdot
b)=\nonumber\\
&=& \sum \varphi (a)(h_{(1)})\varphi (b)(h_{(2)}) =
\varphi(a)\ast \varphi (b) (h), \nonumber
\end{eqnarray}
for all $h\in H$. Therefore $\varphi$ is multiplicative.

For the third claim, let $h,k\in H $ and $a,b \in A$; then
\begin{eqnarray}
\varphi (b(h\cdot a))(k) &=& k\cdot (b(h\cdot a))=\sum (k_{(1)} \cdot b)
(k_{(2)}\cdot(h\cdot a))=\nonumber\\
&=& \sum (k_{(1)} \cdot b)(k_{(2)} \cdot \um)
(k_{(3)}h\cdot a)=\nonumber\\
&=& \sum (k_{(1)} \cdot b)
(k_{(2)}h\cdot a)=\nonumber\\
&=& \sum \varphi (b) (k_{(1)}) \varphi (a) (k_{(2)} h) = \nonumber\\
& = & 
\sum \varphi (b) (k_{(1)}) (h\triangleright \varphi (a))
(k_{(2)})=\nonumber\\
&=& \varphi (b) \ast (h\triangleright \varphi (a)) (k), \nonumber
\end{eqnarray}
$\forall k\in H$. Therefore, $\varphi (h\cdot a)= \varphi (\um)
(h\triangleright \varphi(a))$. The second item follows from this one putting $b=\um$.
\end{proof}

This result suggests that the partial action on $A$ is equivalent to an induced action on $\varphi(A)$, but $\varphi(A)$ must also be a right ideal of an $H$-module algebra; while this may not hold in $Hom_k(H,A)$, it will be true in a certain subalgebra.

\begin{lemma} \label{produto}Let $B$ be  an $H$-module algebra, $x,y \in B$ and $h,k \in H$. Then\\ (i) $(h \hits x)y = \sum h_{(1)} \hits (x (S(h_{(2)})\hits y))$\\
(ii) $(h \hits x)(k \hits y) = \sum h_{(1)} \hits (x (S(h_{(2)})k) \hits y))$.
\end{lemma}
 
\noindent A proof of (i) can be found in (\cite{romenos},lemma 6.1.3) and (ii) is a straightforward consequence of (i).

\begin{prop}\label{ideal}
Let $\varphi: A \vai Hom_k(H,A)$ be as above and consider the $H$-submodule $B = H \hits \varphi(A)$. \\
(i) $B $ is an $H$-module subalgebra of $Hom_k(H,A)$ .\\
(ii) $\varphi(A)$ is a right ideal in $B $ with unity $\varphi(\um)$.
\end{prop}
\begin{proof}
(i) Clearly, $B$ is a $H$-submodule of $Hom_k(H,A)$. Now, given $h \hits \varphi(a)$ and $k \hits \varphi(b) \in H \hits \varphi(A)$, we have 
\begin{eqnarray}
(h \hits \varphi(a))( k \hits \varphi(b)) &  = &  \sum h_{(1)} \hits (\varphi(a) (S(h_{(2)})k) \hits \varphi(b))) \nonumber \\
& = & \sum h_{(1)} \hits \varphi(a (S(h_{(2)})k) \cdot b)) \nonumber 
\end{eqnarray}
and this shows that $B$ is also a  subalgebra.

(ii) This follows by lemma \ref{tecnico}, since $ \varphi (b) \ast (h\triangleright \varphi(a)) = \varphi (b (h\cdot a))$. 
\end{proof}

Lemmas \ref{tecnico}, \ref{produto} and proposition \ref{ideal} prove the existence of enveloping actions. 
\begin{teo}\label{existence}
Let $A$ be a partial $H$-module algebra and let $\varphi: A \vai Hom_k(H,A)$ be the map given by $\varphi(a)(h) = h \cdot a$, and let $B = H \hits \varphi(A)$; then $(B, \varphi)$ is an enveloping action of $A$.
\end{teo}

We will call $(B,\varphi)$ the {\it standard} enveloping action of $A$.

A special case which will be useful for further results is the case when
$\varphi (A)$ is a bilateral ideal of $B$. When this occurs, the element $\varphi
(\um)$ becomes automatically a central idempotent in $B$ and we have also the
following result:

\begin{prop} \label{bilateral} Let $A$ be a partial $H$ module algebra and let $\varphi: A \vai Hom_k(H,A)$  and $B= H \hits \varphi(A)$ be as before. Then  $\varphi (A)\trianglelefteq B$ if and only if
\[
h\cdot (k\cdot a)=\sum (h_{(1)} k\cdot a)(h_{(2)} \cdot \um), \qquad \forall
a\in A, \quad \forall h,k\in H.
\]
\end{prop}

\begin{proof} Suppose that $\varphi(A)$ is an ideal of $B$. We already know that $\forall k\in H$ and $\forall a\in A$ we
  have
\[
\varphi (k\cdot a)= \varphi (\um) \ast (k\triangleright \varphi(a))=
(k\triangleright \varphi(a))\ast\varphi (\um)
\]
Then, these two functions coincide for all $h\in H$:
\[
\varphi (k\cdot a) (h)=  (k\triangleright \varphi(a))\ast \varphi (\um)(h).
\]
The left hand side of the previous equality leads to
\begin{equation}
\label{lefthandside}
\varphi (k\cdot a) (h)=h\cdot (k \cdot a).
\end{equation}
While the right hand side gives
\begin{eqnarray}
\label{righthandside}
 (k\triangleright \varphi(a))\ast \varphi (\um)(h) &=&
\sum (k\triangleright \varphi(a)) (h_{(1)})\ast
\varphi (\um) (h_{(2)}) =\nonumber\\
&=& \sum \varphi(a)(h_{(1)} k) \varphi (\um)(h_{(2)})=\nonumber\\
&=&\sum (h_{(1)} k\cdot a)(h_{(2)} \cdot \um).
\end{eqnarray}
Combining the expressions (\ref{lefthandside}) with (\ref{righthandside}), we
have the result.

Conversely, suppose that $h\cdot (k\cdot a)=\sum (h_{(1)} k\cdot a)(h_{(2)} \cdot \um)$ holds for all $a \in A$ and $h,k \in H$. Equations (\ref{lefthandside}) and (\ref{righthandside}) show that 
\[
\varphi(\um) (k \hits \varphi(a)) = (k \hits \varphi(a)) \varphi(\um) 
\] 
for every $a \in A$ and $k \in H$, i.e., $\varphi(\um)$ is a central idempotent in $B$; therefore $\varphi(A) = \varphi(\um) B$ is an ideal in $B$.
\end{proof}

In \cite{dok} the authors proved the uniqueness of the enveloping action for
a partial action of a group on a unital algebra $A$. In this case, we
have seen that the existence of an enveloping action depends on the fact that
every ideal $D_g$ is endowed with an unity $1_g$. The idea to prove the 
uniqueness is to suppose that there exist two algebras $B$ and $B'$ with actions
$\beta$ and $\beta '$ of the group $G$, respectively, and embeddings 
$\varphi : A\rightarrow B$ and $\varphi ':A\rightarrow B'$ such that 
the partial action of $G$ on $A$ is equivalent to an admissible restriction of
both $\beta $ and $\beta '$. Then one defines a map $\phi : B'\rightarrow B$ by 
$\phi (\beta_g ' (\varphi ' (a)))=\beta_g (\varphi (a))$. The main difficulty
in this theorem is to prove that this map $\phi$ is well defined as a linear
map. This is achieved by two results: first, for each $g\in G$ the
subspace $\beta_g (\varphi (A))$ is an ideal with unity in $B$ (the same
occuring in $B'$), and second, the sum of a finite number of ideals with unity
is also an ideal with unity.

In the Hopf algebra context, the situation is a bit different. Let $H$ be a Hopf
algebra acting partially on a unital algebra $A$ and $(B,\theta)$ be an
enveloping action. By definition of enveloping action, $\theta (A)$ must be a right ideal of $B$
but, unless $h\in H$ is a grouplike element, it is no longer true that the subspaces $h\triangleright \theta (A)$ are
right ideals of $B$; neither it is
true that the elements $h\hits \theta (\um)$ behave as some kind of unity.

In fact, a simple example shows that one doesn't have uniqueness of 
the enveloping action unless stronger conditions have been assumed. Let us
recall the partial action obtained from the adjoint action of the Sweedler
Hopf algebra $H_4$ on itself. The partial action is constructed when we reduce
to the quotient algebra $\overline{B} =H_4 /\langle h_2 \rangle$ and taking
the residual action restricted to the right ideal 
$A= \langle e_1 \rangle$. One enveloping action is given by $(B,i)$, where  $B = H_4 \hits A$ and $i: A \vai B$ is the inclusion. Note that 
$H_4 \hits A = \langle \overline{e_1}, \overline{h_1} \rangle$, 
which is isomorphic to the algebra 
\[
\left[
\begin{array}{cc} 
k & 0 \\
k & 0 
\end{array} 
\right] .
\]
Hence, $B = \langle \overline{e_1}, \overline{h_1} \rangle$ is an enveloping 
algebra of $A$. Nevertheless, the induced action 
on $A$ is total, as we have seen, and hence $(A,\mbox{Id}_A)$ is a ``smaller'' enveloping action. 

In order to clarify these matters we are going to relate each enveloping action of $A$  with the standard enveloping action $(B,\varphi)$ given in theorem \ref{existence}. Suppose that $(B',\theta)$ is an enveloping action for $A$. Define the map 
$\Phi :B' \rightarrow Hom_k (H,A)$ by 
\[
\Phi (\sum_{i=1}^n h_i \triangleright \theta (a_i ))=\sum_{i=1}^n h_i 
\triangleright \varphi (a_i),
\]
where $\varphi :H\rightarrow A$ and the $H$-action on $Hom_k(H,A)$ are the same as before.

\begin{teo} \label{morfismo} If $(B', \theta)$ is an enveloping action of $A$, then
  the map $\Phi$ is an $H$-module algebra morphism onto 
$B = H\triangleright \varphi (A)$.
\end{teo}

\begin{proof} First we have to check that $\Phi$ is well defined as linear
  map. In order to do this, it is enough to prove that if $\sum_{i=1}^n h_i \triangleright \theta (a_i )=0 $ then $\sum_{i=1}^n h_i \triangleright \varphi (a_i )=0 $. So, suppose that 
$x=\sum_{i=1}^n h_i \triangleright \theta (a_i )=0 $; then,  for all $k \in H$,
\[
0 = \theta(\um) (k \hits x) = \theta(\um) \sum_{i=1}^n kh_i \triangleright \theta
(a_i )=\sum_{i=1}^n kh_i \cdot \theta(a_i) = \theta( \sum_{i=1}^n kh_i \cdot \theta(a_i)) 
\]
and, since $\theta$ is injective, it follows that 
\[
\sum_{i=1}^n kh_i \cdot a_i =0
\]
for all $k \in H$. Hence 
\[
\Phi (x) (k) = \sum_{i=1}^n h_i \triangleright \varphi (a_i )(k) = \sum_{i=1}^n  \varphi (a_i )(kh_i ) = \sum_{i=1}^n kh_i \cdot a_i =0
\]for all $k \in H$, which means that $\Phi (x)=0$ and that $\Phi$ is well defined.

$\Phi$ is linear by construction;  we have to show that it is an algebra
morphism. Given $h,k \in H$ and $a,b \in A$,
\begin{eqnarray}
\Phi ((h\triangleright \theta (a))(k\triangleright \theta (b)))  &=&
\Phi (\sum h_{(1)} \triangleright 
(\theta (a)(S(h_{(2)})k\triangleright \theta (b))))=\nonumber\\
&=& \Phi( \sum h_{(1)} \triangleright 
\theta (a(S(h_{(2)})k\cdot b)))=\nonumber\\
&=& \sum h_{(1)} \triangleright \varphi (a(S(h_{(2)})k\cdot b))=\nonumber\\
&=& \sum h_{(1)} \triangleright 
(\varphi (a)(S(h_{(2)})k\triangleright \varphi (b)))=\nonumber\\
&=&  \sum (h_{(1)} \triangleright\varphi (a)) \ast
(h_{(2)} \triangleright(S(h_{(3)})k\triangleright \varphi (b)))=\nonumber\\
&=&  \sum (h_{(1)} \triangleright\varphi (a)) \ast
(h_{(2)} S(h_{(3)})k\triangleright \varphi (b))=\nonumber\\
&=& (h\triangleright \varphi (a)) \ast (k\triangleright \varphi (b))=\nonumber\\
&=& \Phi (h\triangleright \theta (a)) \ast \Phi (k\triangleright \theta (b)).\nonumber
\end{eqnarray}

The fact that $\Phi$ is a morphism of left $H$ modules is easily obtained by
the definition of this map. And since
\[
\sum_{i=1}^n h_i \triangleright \varphi (a_i)= 
\Phi (\sum_{i=1}^n h_i \triangleright \theta (a_i )).
\]
$Phi$ is a surjective map from $B'$ onto $B = H \hits \varphi(A)$.
\end{proof}

Now, {\it injectivity} of $\Phi$ will follow only if whenever $\sum_{i=1}^n kh_i \cdot a_i =0$ for all $k \in H$, then $\sum_{i=1}^n h_i \hits \theta(a_i)$. This motivates the following definition. 

\begin{defi}  Let  $A$ be a partial $H$-module algebra. An enveloping action $(B,\theta)$ of $A$ is minimal if, for every $H$-submodule $M$ of $B$,  
$\theta(\um ) M = 0$ implies $M =0$.
\end{defi}

This concept does not appear in the theory of  partial group actions because, in this case, {\it every enveloping action in minimal}. In fact, consider $H =kG$, $B$ a $H$-module algebra, $A$ a right ideal with unity in $B$ such that $B = kG \hits A$. Suppose that $\sum_i g_i \hits a_i \in B$ satisfies
\[
\sum_i hg_i \hits a_i =0, \forall h \in G.
\] 

Let $(h \hits b) \in B$, with $h \in G$. Then, by the lemma \ref{produto}
\[
(h \hits b) (\sum_i g_i \hits a_i) =  
h \hits (b [\sum_i (h^{-1 }g_i \cdot a_i)]) =  0
\] 

By the results of Exel and Dokuchaev \cite{dok}, the ideal 
$(g_1 \hits A) + \cdots + (g_n \hits A)$ has a unity; since 
$B (\sum_i g_i \hits a_i) =0$, we conclude that $\sum_i g_i \hits a_i =0$. 

\bigskip

On the other hand, in the case of the action of $H_4$ on $\overline{B}=H_4/\langle h_2\rangle$, we had $A = \langle \ee_1 \rangle$, $B = H_4 \hits A = \langle \ee_1, \hh_1 \rangle$, and $\hh_1 = h_1 \hits \ee_1$; then $\Phi(B) = \varphi(A)$, because $\Phi(\hh_1) = h_1 \hits \varphi(e_1)$ is the zero function. This shows that an enveloping action does not need to be minimal. However, regardless of $A$ or $H$, the standard enveloping action is always minimal.

\begin{lemma}\label{minimal}
Let $\varphi: A \vai Hom_k(H,A)$ be as above and consider the $H$-submodule $B = H \hits \varphi(A)$. Then, $(B,\varphi)$ is a minimal enveloping action of $A$.
\end{lemma}
\begin{proof}
It is enough to check that the minimality  condition holds for cyclic submodules. Let  \mbox{$M = H \hits (\sum_{i=1}^n h_i \hits \varphi(a_i))$} and suppose  that $\varphi(\um) \ast M =0$. This means that 
\[0 = \varphi(\um) \ast (\sum_{i=1}^n kh_i \hits \varphi(a_i))  = \sum_{i=1}^n \varphi(kh_i \cdot a_i) = \varphi(\sum_{i=1}^n kh_i \cdot a_i)\] 
for each $k \in H$. Since $\varphi$ is injective, $\sum_{i=1}^n kh_i \cdot a_i = 0$ for every $k \in H$. But then
\[
\sum_{i=1}^n h_i \hits \varphi(a_i)(k) = \sum_{i=1}^n \varphi(a_i)(kh_i) = \sum_{i=1}^n \varphi(kh_i \cdot a_i) = 0
\]
for each $k \in H$, and we conclude that $\sum_{i=1}^n h_i \hits \varphi(a_i) =0$.
\end{proof}

By Theorem \ref{morfismo} and lemma \ref{minimal} we conclude: 
\begin{teo}
Every partial $H$-module algebra has a minimal enveloping action, and any two minimal enveloping actions of $A$ are isomorphic as $H$-module algebras. Moreover, if $(B',\theta)$ is an enveloping action, then there is a morphism of $H$-module algebras of $B'$ onto a minimal enveloping action.
\end{teo}
\begin{proof}
Since the standard enveloping action $(B,\varphi)$ is a minimal action and \linebreak $\Phi: B' \vai Hom_k(H,A)$ is a $H$-module algebra isomorphism onto $B$, we just have to prove that it is injective if $(B',\theta) $ is minimal. If 
$\Phi (\sum_{i=1}^n h_i \triangleright \theta (a_i )) =0$, then 
\[
\sum_{i=1}^n h_i \triangleright \varphi (a_i) (k) = \sum_{i=1}^n kh_i \cdot a_i =0
\]
for each $k \in H$, and hence $0 = \theta(\sum_{i=1}^n kh_i \cdot a_i) = \theta(\um) [k \hits (\sum_{i=1}^n h_i \hits \theta(a_i))] $ for every $k \in H$; hence, $\theta(\um)M=0$ for $M = H \hits (\sum_{i=1}^n h_i \hits \theta(a_i)) =0$, implying that $M=0$ and hence that  $\sum_{i=1}^n h_i \hits \theta(a_i) =0$. 
\end{proof}

\section{A Morita context}

In the reference \cite{dok} the authors showed that if a partial group action
$\alpha$ of $G$ on a unital algebra $A$ admits an enveloping action $\beta$
of the same group on an algebra $B$, then the partial crossed product 
$A\rtimes_\alpha G$ is Morita equivalent to the crossed product 
$B\rtimes_\beta G$. They proved this result by constructing a Morita context
between these two crossed products.

\begin{defi} A Morita context is a six-tuple $(R,S,M,N,\tau, \sigma )$ where
\begin{enumerate}[\bf a)]
\item $R$ and $S$ are rings,
\item $M$ is an $R-S$ bimodule, 
\item $N$ is an $S-R$ bimodule, 
\item $\tau :M\otimes_S N\rightarrow R$ is a bimodule morphism, 
\item $\sigma :N\otimes_R M \rightarrow S$ is a bimodule morphism, 
\end{enumerate}
such that
\begin{equation}
\label{assoctau}
\tau (m_1\otimes n) m_2=m_1 \sigma (n\otimes m_2) , \qquad \forall m_1, m_2 \in M, \quad
\forall n \in N , 
\end{equation}
and 
\begin{equation}
\label{assocsigma}
\sigma (n_1 \otimes m) n_2 = n_1 \tau (m, n_2 ), \qquad \forall m\in M, \quad
\forall n_1, n_2 \in N.
\end{equation}
\end{defi}

By a fundamental theorem due to Morita (see, for example \cite{jacobson} on
pages 167-170), if the morphisms $\tau$ and $\sigma$ are surjective, then the
categories ${}_R \mbox{Mod}$ and ${}_S \mbox{Mod}$ are equivalent. In this
case, we say that $R$ and $S$ are Morita equivalent.

In the Hopf algebra case, we shall see that when a Hopf
algebra $H$ acts partially on a unital algebra $A$ and the enveloping action
$(B, \theta)$ is such that $\theta (A)\trianglelefteq B$, then the partial
smash product ${\underline{A\# H}}$ is Morita equivalent to the smash product
$B\# H$. Just remembering \cite{caen06}, the partial smash product is defined
as follows; put an algebra structure in $A\otimes H$ with the product
\[
(a\otimes h)(b\otimes k)=\sum a(h_{(1)}\cdot b)\otimes h_{(2)} k.
\]
The partial smash product is the unital subalgebra of $A\otimes H$ given by
\[
{\underline{A\# H}}=(A\otimes H)(\um \otimes 1_H), 
\]
that is, the subalgebra generated by typical elements of the form
\[
x=\sum a(h_{(1)}\cdot \um)\otimes h_{(2)}, \qquad \forall a\in A, \quad
\forall h\in H.
\]
Although the definition is quite different from that of  the partial crossed product  $A \rtimes_{\alpha} G$, it can be proved that $A \rtimes_{\alpha} G$ is isomorphic to $\underline{A \# kG}$.

\bigskip

Our aim in this section is to construct a Morita context between the partial
smash product ${\underline{A\# H}}$ and the smash product $B\# H$, where $B$
is an enveloping action for the partial left action $\cdot :H\otimes A
\rightarrow A$. For this purpose, we will embed the partial smash product into
the smash product.

\begin{lemma} If a Hopf algebra $H$ acts partially on a unital algebra $A$
  and $(B, \theta)$ is an enveloping action, then there is an algebra
  monomorphism between the partial smash product ${\underline{A\# H}}$ and the
  smash product $B\# H$.
\end{lemma}

\begin{proof} Define first a map 
\[
\begin{array}{rccc}
\widetilde{\Phi} : & A\times H & \rightarrow & B\# H.\\
\, & (a,h) & \mapsto & \theta (a) \# h
\end{array}
\]
It is easy to see that $\widetilde{\Phi}$ is a bilinear map, then, by the
universal property of the tensor product $A\otimes H$ we define a $k$-linear
map
\[
\begin{array}{rccc}
\Phi : & A\otimes H & \rightarrow & B\# H.\\
\, & a\otimes h & \mapsto & \theta (a) \# h
\end{array}
\]

Now, let us check that $\Phi$ is a morphism of algebras
\begin{eqnarray}
\Phi ((a\otimes h)(b\otimes k) &=& \Phi (\sum a(h_{(1)} \cdot b) \otimes
h_{(2)}k) =\nonumber\\
&=& \sum \theta (a(h_{(1)} \cdot b))\# h_{(2)}k=\nonumber\\
&=& \sum \theta (a) (h_{(1)} \triangleright \theta(b)) \# h_{(2)}
k=\nonumber\\
&=& (\theta (a) \# h)(\theta (b) \# k)= \Phi (a\otimes h) \Phi (b\otimes
k). \nonumber
\end{eqnarray}

Next, we must verify that $\Phi$ is injective. For this, take 
\[ 
x=\sum_{i=1}^n a_i \otimes h_i \in \mbox{ker}\Phi , 
\]
and, without loss of generality, choose $\{ a_i \}_{i=1}^n$ to be linearly
independent. As $\theta$ is injective, we conclude that the elements
$\theta (a_i ) \in B$ are linearly independent. For each $f\in H^*$, we have
\[
\sum_{i=1}^n \theta (a_i )f(h_i)=0 ,
\]
and then, by the linear independence of $\{ \theta (a_i )\}$, we get 
$f(h_i )=0$, $\forall f\in H^*$. Therefore $h_i =0$ for $i=1, \ldots n$, which
implies that $x=0$ and that $\Phi$ is injective.

As the partial smash product ${\underline{A\# H}}$ is a sub-algebra of
$A\otimes H$, then, it is injectively mapped into $B\# H$ by $\Phi$. A
typical element of the image of the partial smash product is
\begin{eqnarray}
\Phi ((a\otimes h)(\um \otimes 1_H)) &=& \Phi (a\otimes h)\Phi(\um \otimes
1_H) =\nonumber\\
&=& (\theta (a) \# h) (\theta (\um) \# 1_H) =\nonumber\\
&=& \sum \theta (a) (h_{(1)} \triangleright \theta (\um))\#
h_{(2)}.\nonumber
\end{eqnarray}
\end{proof}

Take $M = \Phi (A \otimes H) = \{\sum_{i=1}^n \theta (a_i) \# h_i; a_i \in A, n \in \N \}$; and take
$N$ as the subspace of $B\# H$ generated by the elements 
$\sum (h_{(1)} \hits \theta(a)) \# h_{(2)}$ with $h \in H$, and $a \in A$. Another way to characterize $N$ is as the subspace 
$N=(\um \otimes H)\Phi (A \otimes 1_H )$.

\begin{prop}
Let $H$ be a Hopf algebra with invertible antipode, $A$ a partial $H$-module algebra, and suppose that $\theta (A)$ is an ideal of $B$; then $M$ is 
a right $B\# H$ module and $N$ is a left $B \# H$ module.
\end{prop}
\begin{proof}
In order to see that $M$ is a right $B\# H$ module, let $\theta (a) \#h \in M$ and $b \# k \in B \# H$. Then 
\[
(\theta (a) \# h)( b \# k) =  \sum \theta(a) (h_{(1)} \hits b) \# h_{(2)} k,
\]
which lies in $\Phi (A \otimes H)$ because $\theta (A)$ is a right ideal in
$B$. 

Now, to prove that $N$ is a left $B\# H$ module, let 
$\sum (h_1 \hits \theta(a)) \# h_2 $ be a generator of $N$.
\begin{eqnarray}
& \, & ( b \# k)(\sum (h_{(1)} \hits \theta(a)) \# h_{(2)})   = \nonumber\\ 
&=& \sum b (k_{(1)}h_{(1)} \hits \theta(a)) \# k_{(2)}h_{(2)} = \nonumber \\
 & = & \sum b (\epsilon(k_{(1)}h_{(1)})k_{(2)}h_{(2)} \hits \theta(a)) 
\# k_{(3)}h_{(3)} =\nonumber\\
 & = & \sum (\epsilon(k_{(1)}h_{(1)}) 1_H \hits b) 
(k_{(2)}h_{(2)} \hits \theta(a)) \# k_{(3)}h_{(3)} =\nonumber\\
 & = & \sum ((k_{(2)}h_{(2)}S^{-1}(k_{(1)}h_{(1)})) \hits b) 
(k_{(3)}h_{(3)} \hits \theta(a)) \# k_{(4)}h_{(4)} =\nonumber\\
 & = & \sum (k_{(2)}h_{(2)}\triangleright(S^{-1}(k_{(1)}h_{(1)}) \hits b))
(k_{(3)}h_{(3)} \hits \theta(a)) \# k_{(4)}h_{(4)} =\nonumber\\
 & = & \sum (k_{(2)}h_{(2)} \hits ((S^{-1}(k_{(1)}h_{(1)}) \hits
 b)\theta(a))) \# k_{(3)}h_{(3)}.
\nonumber
\end{eqnarray}
Each term 
$(S^{-1}(k_{(1)}h_{(1)}) \hits b) \theta(a)$ is in $\theta(A)$ because
$\theta(A)$ is an ideal of $B$. It follows that   $N$ is  a
left $B\# H$ ideal.
\end{proof}

We can define a left ${\underline{A\# H}}$ module structure on $M$ and
 a right ${\underline{A\# H}}$ module structure on $N$ induced by the
monomorphism $\Phi$, that is,
\[
(\sum a(h_{(1)}\cdot \um) \otimes h_{(2)})\blacktriangleright (\theta (b) \# k)
=(\sum \theta (a)(h_{(1)}\triangleright \theta (\um)) \# h_{(2)})
(\theta (b) \# k), 
\]
and
\begin{eqnarray}
&\, & (\sum k_{(1)} \triangleright \theta (b) \# k_{(2)})\blacktriangleleft 
(\sum a(h_{(1)}\cdot \um) \otimes h_{(2)})=\nonumber\\ 
&=& (\sum k_{(1)} \triangleright
\theta (b) \# k_{(2)})
(\sum \theta (a)(h_{(1)}\triangleright \theta (\um)) \# h_{(2)}).\nonumber
\end{eqnarray}

\begin{prop} Under the same hipotheses of the previous proposition, $M$ is
  indeed a 
left ${\underline{A\# H}}$ module with the map $\blacktriangleright$ and 
$N$ is indeed a right ${\underline{A\# H}}$ module with the map 
$\blacktriangleleft$.
\end{prop}

\begin{proof} Let us first verify that 
$ {\underline{A \# H}} \blacktriangleright  M \subseteq M$. 
\begin{eqnarray}
&\, & (\sum a(h_{(1)}\cdot \um) \otimes h_{(2)})
\blacktriangleright (\theta (b) \# k)=\nonumber\\
&=& (\sum \theta (a)(h_{(1)}\triangleright \theta (\um)) \# h_{(2)})
(\theta (b) \# k) =\nonumber\\
&=& \sum \theta (a) (h_{(1)} \cdot \theta (\um)) 
(h_{(2)} \hits \theta (b))\# h_{(3)} k =\nonumber\\
& = & \sum \theta (a)(h_{(1)} \cdot \theta (\um) \theta (b)) \# h_{(2)} k= \nonumber\\
& = & \sum \theta (a)(h_{(1)} \cdot \theta (b)) \# h_{(2)} k, \nonumber
\end{eqnarray}
which lies inside $M$ because $\theta (A)$ is a right ideal of $B$.

Now, we have to verify that 
$N \blacktriangleleft  {\underline{A \# H}}\subseteq N$. 
\begin{eqnarray}
& \, & (\sum k_{(1)} \triangleright \theta (b) \# k_{(2)})\blacktriangleleft 
(\sum a(h_{(1)}\cdot \um) \otimes h_{(2)}) =\nonumber\\ 
&=& (\sum k_{(1)} \triangleright \theta (b) \# k_{(2)})
(\sum \theta (a)(h_{(1)}\triangleright \theta (\um)) \# h_{(2)})=\nonumber\\
&=& \sum (k_{(1)} \triangleright \theta (b))
(k_{(2)} \triangleright (\theta (a)(h_{(1)}\triangleright \theta (\um))) 
\# k_{(3)} h_{(2)} =\nonumber\\
&=& \sum (k_{(1)} \triangleright \theta (b))
(k_{(2)} \triangleright \theta (a))
(k_{(3)}\triangleright(h_{(1)}\triangleright \theta (\um))) 
\# k_{(4)} h_{(2)} =\nonumber\\
&=& \sum (k_{(1)} \triangleright \theta (ba))
(k_{(2)}h_{(1)}\triangleright \theta (\um)) 
\# k_{(3)} h_{(2)} =\nonumber\\
&=& \sum (k_{(1)}\epsilon (h_{(1)}) \triangleright \theta (ba))
(k_{(2)}h_{(2)}\triangleright \theta (\um)) 
\# k_{(3)} h_{(3)} =\nonumber\\
&=& \sum (k_{(1)} h_{(2)} S^{-1} (h_{(1)}) \triangleright \theta (ba))
(k_{(2)}h_{(3)}\triangleright \theta (\um)) 
\# k_{(3)} h_{(4)} =\nonumber\\
&=& \sum (k_{(1)} h_{(2)}\triangleright( S^{-1} (h_{(1)}) \triangleright \theta (ba)))
(k_{(2)}h_{(3)}\triangleright \theta (\um)) 
\# k_{(3)} h_{(4)} =\nonumber\\
&=& \sum k_{(1)} h_{(2)}\triangleright(( S^{-1} (h_{(1)}) \triangleright \theta (ba))\theta (\um)) 
\# k_{(2)} h_{(3)} \nonumber\\
& = & \sum k_{(1)} h_{(2)}\triangleright\theta( S^{-1} (h_{(1)}) \cdot (ba))
\# k_{(2)} h_{(3)} \nonumber
\end{eqnarray}
where in the last passage we used the fact that   
\[
\theta(t \cdot x) = \theta(\um) (t \hits \theta(x)) = (t \hits \theta(x)) \theta(\um).
\]
which holds because $\theta(\um)$ is a central idempotent.
\end{proof}

The last ingredient for a Morita context is the definition of two bimodule
morphisms
\[
\tau : M\otimes_{B\# H} N \rightarrow {\underline{A\# H}} \cong \Phi
({\underline{A\# H}})\subseteq B\# H
\]
and
\[
\sigma : N\otimes_{\underline{A\# H}} M \rightarrow B\# H .
\]
As $M$, $N$ and ${\underline{A\# H}}$ are viewed as subalgebras of $B\# H$,
these two maps can be taken as the usual multiplication on $B\# H$. The
associativity of the product assures us that these maps are bimodule morphisms and
satisfy the associativity conditions (\ref{assoctau}) and
(\ref{assocsigma}). The following proposition shows us that the maps $\tau$
and $\sigma$ are indeed well defined, and furthermore, that they are
surjective, proving the Morita equivalence between these two smash products. 
 
\begin{prop}
$MN = \Phi (\underline{A \# H})$, $NM = B \# H$
\end{prop}

\begin{proof}
Let's see first that  $MN \subseteq \Phi (\underline{A \# H})$. 

\begin{eqnarray}
& \, & (\theta (a) \# h) (\sum (k_{(1)} \hits \theta (b)) \# k_{(2)})
=\nonumber\\
& = & \sum \theta (a) h_{(1)} \hits (k_{(1)} \hits \theta (b)) \# h_{(2)} k_{(2)} =\nonumber\\
& = & \sum \theta (a) (h_{(1)} k_{(1)} \hits \theta (b)) \# h_{(2)} k_{(2)} =\nonumber\\
& = & \sum \theta (a)  (h_{(1)} k_{(1)} \hits \theta (b)) 
(h_{(2)} k_{(2)} \hits \theta(\um))\# h_{(3)} k_{(3)} =\nonumber\\
& = & \sum \theta (a  (h_{(1)} k_{(1)} \cdot b)) 
((h_{(2)} k_{(2)})_{(1)} \hits \theta (\um))\# (h_{(2)} k_{(2)})_{(2)}. \nonumber
\end{eqnarray}

which is an element of $\Phi ({\underline{A\# H}})$. 

Since 
\[
\sum \theta (a) (h_{(1)} \triangleright \theta (\um) )\# h_{(2)} =(\theta
(a)\# h)(\theta (\um) \# 1_H)
\]
and $\theta (\um) \# 1_H \in N$, it follows that $MN = \Phi ( \underline{A \# H})$.

In order to prove $NM = B \# H$ we just have to show that every element of the form $(h \hits \theta (a)) \# k$ is in $NM$, because this is a generating set for $B \# H$ as a vector space.  We claim that
\[
(h\hits \theta (a))\# k =
\sum ((h_{(1)} \hits \theta (a)) \# h_{(2)})
(\theta (\um) \# S(h_{(3)}) k).
\]
This can be easily seen as follows:
\begin{eqnarray}
& \, & \sum ((h_{(1)} \hits \theta (a)) \# h_{(2)})
(\theta (\um) \# S(h_{(3)}) k =\nonumber\\ 
&=& \sum ((h_{(1)} \hits \theta (a)) 
(h_{(2)} \hits \theta (\um)) \# h_{(3)}S(h_{(4)}) k=\nonumber\\
&=& \sum ((h_{(1)} \hits \theta (a)\theta (\um)) 
 \# \epsilon (h_{(2)}) k=\nonumber\\
&=&  ((\sum h_{(1)} \epsilon (h_{(2)}))\hits \theta (a)) 
 \# k=\nonumber\\
&=& (h\hits \theta (a))\# k. \nonumber
\end{eqnarray}

Therefore $NM=B\# H$.
\end{proof}

In conclusion, we have constructed a Morita context for the two smash products and
proved that this Morita context gives us a Morita equivalence between these
two algebras.

\section{Enveloping coactions}

Following \cite{caen06}, we define a partial right coaction of a Hopf 
algebra $H$ on a algebra $A$ to be a linear map $\prho: A \vai A \otimes H$ 
such that
\begin{eqnarray}
\label{coaction}
\mbox{1) }&& \prho (ab) =\prho (a) \prho (b), \; \forall a,b\in A, \nonumber\\ 
\mbox{2) }&& (I\otimes \epsilon )\prho (a)=a, \; \forall a\in A, \nonumber\\
\mbox{3) }&& (\prho \otimes I)\prho (a) =(\prho (\um) \otimes 1_H )
((I\otimes \Delta)\prho (a)), \; \forall a\in A. 
\end{eqnarray}
We will denote 
\[
\prho (a) =\sum a^{({\underline{0}})} \otimes a^{({\underline{1}})}.
\]
Note that, in this notation, we can rewrite the items 1), 2) and 3) of
(\ref{coaction}) as
\begin{eqnarray}
\mbox{1) }&& \sum (ab)^{({\underline{0}})} \otimes (ab)^{({\underline{1}})} =
a^{({\underline{0}})}b^{({\underline{0}})} \otimes
a^{({\underline{1}})}b^{({\underline{1}})} , \nonumber\\
\mbox{2) }&& \sum a^{({\underline{0}})} \epsilon(a^{({\underline{1}})}) =a,
\nonumber\\
\mbox{3) }&& \sum a^{({\underline{0}})({\underline{0}})}\otimes 
a^{({\underline{0}})({\underline{1}})} \otimes a^{({\underline{1}})} =
\sum \um ^{({\underline{0}})}a^{({\underline{0}})} \otimes 
\um ^{({\underline{1}})}{a^{({\underline{1}})}}_{(1)} \otimes 
{a^{({\underline{1}})}}_{(2)}. \nonumber
\end{eqnarray}

The simplest example of a partial coaction can be given as a restriction of a
coaction of $H$ on a right $H$ comodule algebra $B$.

\begin{prop}
Let $B$ be a right $H$-comodule algebra with coaction 
$\rho : B\rightarrow B\otimes H$. If $A$ is a right ideal of $B$ with 
a unity $\um $, then the map
\[
\prho(a) = (\um  \otimes 1_H) \rho (a) = 
\sum \um  a^{(0)} \otimes a^{(1)}
\] 
defines a
partial right coaction of $H$ on $A$.
\end{prop}
\begin{proof} 
For the first item of (\ref{coaction}) we have, for all $a \in A$
\begin{eqnarray}
\prho (ab) &=& (\um  \otimes 1_H)(\rho (ab))=(\um  \otimes 1_H)(\rho
(a)\rho(b))=\nonumber\\
&=& (\sum \um  a^{(0)} \otimes a^{(1)})( \sum b^{(0)} \otimes b^{(1)}) = \nonumber\\ &=&
\sum \um  a^{(0)}b^{(0)} \otimes a^{(1)}b^{(1)}=\nonumber\\
&=& \sum \um  a^{(0)} \um  b^{(0)} \otimes a^{(1)}b^{(1)}= \nonumber\\
&=& (\sum \um  a^{(0)} \otimes a^{(1)})( \sum \um  b^{(0)} 
\otimes b^{(1)})= \nonumber\\
&=& \prho (a) \prho (b)\nonumber
\end{eqnarray}
where we used  $ \um  x= \um  x \um \in A$ in the fourth line.

Item 2) is easily established for every $a\in A$:
\[
(I\otimes \epsilon )\prho (a)=\sum \um  a^{(0)}\epsilon(a^{(1)})=\um 
a=a.
\]

Finally checking item 3), we have 
\begin{eqnarray}
(\prho \otimes I) \prho (a) & = & \sum (\um \otimes 1_H \otimes I)(\rho \otimes I) (\um \otimes 1_H)(a^{(0)} \otimes a^{(1)}) = \nonumber\\
& = & \sum (\um \otimes 1_H \otimes I) ((\um a)^{(0)} \otimes (\um a)^{(1)} \otimes a^{(2)}) = \nonumber\\
& = & \sum (\um \otimes 1_H \otimes I) (\um^{(0)}a^{(0)} \otimes \um^{(1)}a^{(1)} \otimes a^{(2)}) = \nonumber\\
& = &  (\um \otimes 1_H \otimes I) (\rho(\um) \otimes I) (\rho \otimes I) \rho (a) = \nonumber\\
& = & (\um \otimes 1_H \otimes I) (\rho(\um) \otimes I) (I \otimes \Delta) \rho (a) = \nonumber\\
& = & (\prho(\um) \otimes I) (I \otimes \Delta) \rho (a).\nonumber
\end{eqnarray}
Therefore $\prho$ is a partial coaction. 
\end{proof}

In \cite{lomp07} the author proved that if a finite dimensional Hopf algebra $H$
acts partially on a unital algebra $A$, then there exists a partial coaction
of the dual Hopf algebra $H^*$ on $A$ by the coaction
\[
\prho (a) = \sum_{i=1}^{n} (b_i \cdot a) \otimes p_i, 
\]
where $\{ b_i \}_{i=1, \ldots n}$ is a basis for $H$ and 
$\{ p_i \}_{i=1, \ldots n}$ is its dual basis in $H^*$. In fact, one can push
forward this result to a more general context. We recall the definition of a pairing of Hopf algebras.
\begin{defi} A pairing between two Hopf algebras $H_1$ and $H_2$ is a linear map
\[
\begin{array}{rccc}
\langle , \rangle : & H_1 \otimes H_2 & \rightarrow & k\\
\, & h\otimes f & \mapsto & \langle h,f \rangle 
\end{array}
\]
such that:
\begin{enumerate}
\item[(i)] $\langle hk,f \rangle=\langle h\otimes k,\Delta (f) \rangle$.
\item[(ii)] $\langle h,fg \rangle=\langle \Delta (h),f\otimes g \rangle$.
\item[(iii)] $\langle h,1_{H_2} \rangle =\epsilon (h)$.
\item[(iv)] $\langle 1_{H_1},f \rangle =\epsilon (f)$.
\end{enumerate}
A pairing is said to be nondegenerate if the following conditions hold:
\begin{enumerate}
\item If $\langle h,f \rangle=0$, $\forall f\in H_2$ then $h=0$.
\item If $\langle h,f \rangle=0$, $\forall h\in H_1$ then $f=0$. 
\end{enumerate}
\end{defi}
Let $H_1$ and $H_2$ two Hopf
algebras dually paired with a nondegenerate pairing, and suppose also that $H_1$ acts partially on an algebra $A$ in such a way that
for all $a\in A$ the subspace $H_1 \cdot a$ has finite dimension (this is a
requirement that is analogous to the case of rational modules). Then we have
the following result:

\begin{teo}  Let $H_1$ and $H_2$ two Hopf
algebras dually paired with a nondegenerate pairing, and suppose that $H_1$ acts
partially on an algebra $A$ in such a manner that $\mbox{dim}(H_1 \cdot a)<\infty$ for all $a \in A$. Then
$H_2$ coacts partially on $A$ with a coaction $\prho :A\rightarrow A\otimes
H_2$ defined by
\[
(I\otimes h)\prho (a)= h\cdot a, \qquad \forall h\in H_1
\]
where $I\otimes h :A\otimes H_2 \rightarrow A$ is given by 
\[
(I\otimes h)(\sum_{i=1}^{n} a_i \otimes f_i)=\sum_{i=1}^{n} a_i \langle h,f_i
\rangle.
\]
\end{teo}

\begin{proof} The condition that $\mbox{dim}(H_1 \cdot a)<\infty$, $\forall
  a\in A$, implies that given $a\in A$ there exist elements 
$a_1, \ldots a_n\in A$ and $\lambda_1, \ldots, \lambda_n \in H_1^\ast$ such that 
\[
h\cdot a=\sum_{i=1}^{n}  \lambda_i(h)a_i,
\]
where we may choose the elements $a_i$ linearly
independent. The
map $T_a :H\rightarrow A$, given by, $T_a(h) =  h\cdot a$, has a kernel
of codimension $n$, and there is an $n$-dimensional subspace $V_a$ of $H_1$ such that $H_1 = V_a \oplus \ker(T_a)$. Choose a basis $\{h_1, \ldots, h_n\}$ of $V_a$ such that $T_a(h_i) = a_i$ for each $i$;  since the pairing is non-degenerate, there are $f_1, \ldots , f_n \in H_2$ such that $\langle h_i,f_j \rangle = \delta_{i,j}$. It follows that 
for all $h\in H_1$ we have
\[
h\cdot a=\sum_{i=1}^{n} \langle h,f_i \rangle a_i.
\]
Define then $\prho :A \rightarrow A\otimes H_2$ as
\[
\prho (a)=\sum_{i=1}^{n} a_i \otimes f_i.
\]
For each $h\in H_1$,
\[
(I\otimes h)\prho (a) =\sum_{i=1}^{n} a_i \langle h, f_i \rangle =h\cdot a.
\]

We have to verify itens 1), 2) and 3) of (\ref{coaction}). For the first item , let $a,b\in A$ such that 
\[
\prho (a) =\sum_{i=1}^n a_i \otimes f_i, \qquad  \prho (b) =\sum_{j=1}^m b_j \otimes g_j .
\]
Then, for every $h\in H_1$ we have
\begin{eqnarray}
(I\otimes h) \prho (ab) &=& h\cdot (ab) =\sum (h_{(1)} \cdot a)
(h_{(2)} \cdot b) =\nonumber\\
&=& \sum (\sum_{i=1}^n a_i \langle h_{(1)} , f_i \rangle )
(\sum_{j=1}^m b_j \langle h_{(2)} , g_j \rangle )=\nonumber\\
&=& \sum \sum_{i=1}^n \sum_{j=1}^m a_i b_j \langle h_{(1)} , f_i \rangle 
\langle h_{(2)} , g_j \rangle =\nonumber\\
&=&  \sum_{i=1}^n \sum_{j=1}^m a_i b_j \langle h, f_i g_j \rangle =\nonumber\\
&=& (I\otimes h)  \sum_{i=1}^n \sum_{j=1}^m a_i b_j \otimes f_i g_j
=\nonumber\\
&=& (I\otimes h) \prho (a) \prho (b). \nonumber
\end{eqnarray}
As this equality is true $\forall h\in H_1$ and the pairing is nondegenerate,
then we can conclude that
\[
\prho (ab)=\prho (a) \prho (b).
\]

For the second item, we use the equality involving the pairing
\[
\epsilon (f) =\langle 1_{H_1} , f \rangle ,
\] 
then, $\forall a\in A$,
\[
(I\otimes \epsilon)\prho (a)= \sum_{i=1}^n a_i \langle 1_{H_1} , f_i \rangle
=1_{H_1} \cdot a =a.
\]

Finally, for item 3), given any $h,k\in H_1$ we have

\begin{eqnarray}
\label{item3}
&\, &(I\otimes h\otimes k)(\prho \otimes I)\prho (a) = (I\otimes h\otimes k)
\sum_{i=1}^n \prho (a_i ) \otimes f_i = 
\sum_{i=1}^n (I\otimes h)\prho (a_i ) \langle k,f_i \rangle =\nonumber\\
&=&  \sum_{i=1}^n (h\cdot a_i)\langle k,f_i \rangle = 
h\cdot ( \sum_{i=1}^n a_i \langle k,f_i \rangle )= h\cdot (k\cdot a)\nonumber\\
&=& \sum (h_{(1)} \cdot \um)(h_{(2)}k\cdot a)=\nonumber\\
&=& \sum ((I\otimes h_{(1)})\prho (\um))
((I\otimes h_{(2)}k)\prho (a)) =\nonumber\\
&=& \sum ((I\otimes h_{(1)})\prho (\um)) 
\left( \sum_{i=1}^n a_i \langle h_{(2)}k,f_i \rangle \right)=\nonumber\\
&=& \sum ((I\otimes h_{(1)})\prho (\um)) 
\left( \sum_{i=1}^n a_i (\sum \langle h_{(2)},f_{i(1)} \rangle 
\langle k,f_{i(2)} \rangle ) \right)=\nonumber\\
&=& (I\otimes k) (\sum ((I\otimes h_{(1)} \otimes I)(\prho (\um)\otimes 1_{H_1})) 
\left( \sum_{i=1}^n a_i (\sum \langle h_{(2)},f_{i(1)} \rangle \otimes 
f_{i(2)} ) \right).
\end{eqnarray}
If we write $\prho (\um)=\sum_{j=1}^p e_j \otimes \varepsilon_j$ then the last
equality in (\ref{item3}) reads
\begin{eqnarray}
&\, & (I\otimes k) (\sum ((I\otimes h_{(1)} \otimes I)(\prho (\um)\otimes 1_{H_1})) 
\left( \sum_{i=1}^n a_i (\sum \langle h_{(2)},f_{i(1)} \rangle \otimes 
f_{i(2)} ) \right) =\nonumber\\
&=& (I\otimes k) \left(\sum ((I\otimes h_{(1)} \otimes I) 
(\sum_{j=1}^p e_j \otimes \varepsilon_j \otimes  1_{H_1})\right)
\left( \sum_{i=1}^n a_i (\sum \langle h_{(2)},f_{i(1)} \rangle \otimes 
f_{i(2)} ) \right) =\nonumber\\
&=& (I\otimes k) \left(\sum ( 
(\sum_{j=1}^p e_j \langle h_{(1)} ,\varepsilon_j \rangle \otimes  1_{H_1})\right)
\left( \sum_{i=1}^n a_i (\sum \langle h_{(2)},f_{i(1)} \rangle \otimes 
f_{i(2)} ) \right) =\nonumber\\
&=& (I\otimes k) \left( \sum_{j=1}^p \sum_{i=1}^n e_j a_i 
\sum \langle h_{(1)} ,\varepsilon_j \rangle  \langle h_{(2)},f_{i(1)} \rangle 
\otimes f_{i(2)} ) \right) =\nonumber\\
&=& (I\otimes k) \left( \sum_{j=1}^p \sum_{i=1}^n e_j a_i 
\sum \langle h,\varepsilon_j f_{i(1)} \rangle 
\otimes f_{i(2)}  \right) =\nonumber\\
&=& (I\otimes k) (I\otimes h \otimes I) \left( \sum_{j=1}^p \sum_{i=1}^n e_j a_i 
\otimes (\sum \varepsilon_j f_{i(1)}  \otimes f_{i(2)} ) \right) =\nonumber\\
&=& (I\otimes h\otimes k)( \sum_{j=1}^p e_j\otimes  \varepsilon_j\otimes
1_{H_1})( \sum_{i=1}^n a_i \otimes \Delta (f_i))=\nonumber\\
&=& (I\otimes h\otimes k)(\prho (\um) \otimes 1_{H_1})((I\otimes \Delta)\prho
(a)).\nonumber
\end{eqnarray}
As this equality is valid for every $h,k\in H_1$, and because of the
nondegeneracy of the pairing, we conclude that
\[
(\prho \otimes I)\prho (a)=(\prho (\um)\otimes 1_{H_1})((I\Delta)\prho (a)),
\]
$\forall a\in A$. 

Therefore, the map $\prho$ is indeed a partial coaction.
\end{proof}

A special case is when we consider the finite dual of a Hopf algebra $H$,
\[
H^\circ =\{ f\in H^* | \exists I\trianglelefteq H, f(I)=0,
\mbox{dim}H/I<\infty \} .
\]
We say that $H^\circ$ separates points in $H$ if the following condition holds:
\[
f(h)=0, \quad \forall h\in H \Rightarrow f=0.
\]
This condition allows us to define a nondegenerate pairing between $H$ and
$H^\circ$, and therefore, by the previous theorem,  we have the following result:

\begin{cor} Let $H$ be a Hopf algebra which acts partially on a unital algebra $A$ such
  that for each $a\in A$ the subspace $H\cdot a \subseteq A$ is finite
  dimensional. If the finite dual $H^\circ$ of $H$ separates points in $H$,
  then there is a partial coaction of $H^\circ$ on $A$ given by
\[
(I\otimes h) \prho (a) =h\cdot a, \qquad \forall a\in A, \quad \forall h\in H.
\]
\end{cor}

Conversely, suppose that a Hopf algebra $H_1$ coacts partially on an 
unital algebra $A$. Suppose also that there exists a pairing (not necessarilly nondegenerate)
between the Hopf algebras $H_1$ and $H_2$. Then we define a map
\[
\begin{array}{rccc}
{\overline{\cdot}} : & H_2 \times A & \rightarrow & A,\\
\, & (f,a) & \mapsto & \sum a^{({\underline{0}})} 
\langle a^{({\underline{1}})} ,f \rangle
\end{array}
\]
where $\prho (a)= \sum a^{({\underline{0}})} \otimes a^{({\underline{0}})}$ is
the expression of the partial coaction of $H_1$ on $A$. It is easy to see that
this map is bilinear; by the universal property of the tensor product, we
can define a linear map
\[
\begin{array}{rccc}
\cdot : & H_2 \otimes A & \rightarrow & A.\\
\, & f\otimes a & \mapsto & \sum a^{({\underline{0}})} 
\langle a^{({\underline{1}})} ,f \rangle
\end{array}
\]

\begin{prop} The map $\cdot $ defined above is indeed a partial action of
  $H_2$ on $A$.
\end{prop}

\begin{proof} Let $f\in H_2$ and $a,b\in A$, then
\begin{eqnarray}
f\cdot (ab) &=& \sum (ab)^{({\underline{0}})} 
\langle (ab)^{({\underline{1}})} ,f \rangle =\nonumber\\
&=&  \sum a^{({\underline{0}})}b^{({\underline{0}})} 
\langle a^{({\underline{1}})}b^{({\underline{1}})} ,f \rangle =\nonumber\\
&=&  (\sum a^{({\underline{0}})} 
\langle a^{({\underline{1}})} ,f_{(1)} \rangle)( \sum b^{({\underline{0}})} 
\langle b^{({\underline{1}})} ,f_{(2)} \rangle) =\nonumber\\
&=& \sum (f_{(1)}\cdot a)(f_{(2)} \cdot b).\nonumber
\end{eqnarray}

Now, let $a\in A$, we have
\begin{eqnarray}
1_{H_2} \cdot a &=& \sum a^{({\underline{0}})} 
\langle a^{({\underline{1}})} ,1_{H_2} \rangle =\nonumber\\
&=& \sum a^{({\underline{0}})} \epsilon( a^{({\underline{1}})})=a. \nonumber
\end{eqnarray}

Finally, for each $f,g\in H_2$ and $a\in A$, and writing 
$\prho (\um)=\sum 1^{({\underline{0}})} \otimes 1^{({\underline{1}})}$, we have
\begin{eqnarray}
f\cdot (g\cdot a) &=& f \cdot (\sum a^{({\underline{0}})} 
\langle a^{({\underline{1}})} ,g \rangle) =\nonumber\\
&=& \sum a^{({\underline{0}})({\underline{0}})} 
\langle a^{({\underline{0}})({\underline{1}})} ,f \rangle 
\langle a^{({\underline{1}})} ,g \rangle =\nonumber\\
&=& \sum 1^{({\underline{0}})}a^{({\underline{0}})} 
\langle 1^{({\underline{1}})}{a^{({\underline{1}})}}_{(1)} ,f \rangle 
\langle {a^{({\underline{1}})}}_{(2)} ,g \rangle =\nonumber\\
&=& \sum 1^{({\underline{0}})}a^{({\underline{0}})} 
\langle 1^{({\underline{1}})} , f_{(1)} \rangle 
\langle {a^{({\underline{1}})}}_{(1)} ,f_{(2)} \rangle 
\langle {a^{({\underline{1}})}}_{(2)} ,g \rangle =\nonumber\\
&=& (\sum 1^{({\underline{0}})}
\langle 1^{({\underline{1}})} , f_{(1)} \rangle )
(\sum a^{({\underline{0}})} 
\langle {a^{({\underline{1}})}} ,f_{(2)} g \rangle)=\nonumber\\
&=& \sum (f_{(1)}\cdot \um)((f_{(2)} g) \cdot a).\nonumber
\end{eqnarray}

These three properties show that $\cdot $ is in fact a partial action of $H_2$
on $A$.
\end{proof}

A natural question is whether it is possible or not to define an enveloping
coaction for a partial coaction of a Hopf algebra $H$ on a unital algebra
$A$. The basic idea is to use the previous proposition to define a partial 
action of the finite dual $H^\circ$
of the Hopf algebra $H$ on $A$, then take an enveloping action
$(B, \theta)$ of this action, and finally, considering the case when $H^\circ$
separates points, to analyse if 
$\mbox{dim} (H^\circ \triangleright \theta(a)) < \infty$, $\forall a\in A$
(i.e., if  $B$ is a rational left $H^\circ$ module); if this holds, one defines
back a coaction $\rho :B\rightarrow B\otimes H$ by the equation
\[
(I\otimes f)\rho (b) =f\triangleright b, \qquad \forall f\in H^\circ \quad
\forall b\in B.
\]

\begin{prop} Let $H$ be a Hopf algebra such that its finite dual $H^\circ$
  separates points. Suppose that $H$ coacts partially on a unital algebra $A$
  with coaction $\prho$ and an 
enveloping action, $(B, \theta)$, of the partial action of $H^\circ$ on $A$ is a
rational left $H^\circ$ module. Then the map $\theta :A\rightarrow B$
intertwins the partial coaction of $H$ on $A$ and the induced partial coaction
of $H$ on $B$ given by
\[
\widetilde{\rho} (b)=(\theta (\um)\otimes 1_H)\rho (b) , \qquad \forall b\in B.
\]
\end{prop}

\begin{proof} We have to show that 
$(\theta \otimes I) \circ \prho =\widetilde{\rho} \circ \theta$. 
Take $a\in A$ and $f\in H^\circ$, then
\begin{eqnarray}
(I\otimes f) (\theta \otimes  I) \prho (a) &=& 
(I\otimes f)(\sum \theta (a^{({\underline{0}})}) \otimes
a^{({\underline{1}})}) =\nonumber\\
&=& \sum \theta (a^{({\underline{0}})}) \langle a^{({\underline{1}})} ,f
\rangle =\nonumber\\
&=& \theta (\sum a^{({\underline{0}})} \langle a^{({\underline{1}})} ,f
\rangle) =\nonumber\\
&=& \theta (f\cdot a) =f\cdot \theta (a) =\nonumber\\
&=& \theta (\um) (f\triangleright \theta (a))= \theta (\um)((I\otimes f)\rho
(\theta (a)))=\nonumber\\
&=& (I\otimes f)((\theta (\um)\otimes 1_H) \rho (\theta (a)))=
(I\otimes f)\widetilde{\rho} (\theta(a)).\nonumber
\end{eqnarray}
As this identity holds for every $f\in H^\circ$ and $H^\circ $ separates points
in $H$, then $(\theta \otimes  I) \prho (a)=\widetilde{\rho} (\theta(a))$,
$\forall a\in A$. Therefore the map $\theta$ intertwins these two partial
coactions.
\end{proof}

Certainly, the conditions on the existence of enveloping coactions are quite
restrictive but, at least, one class of examples can be given where this occurs,
namely, the finite dimensional Hopf algebras.

\begin{prop} A partial coaction of a finite dimensional Hopf algebra $H$ on an
  unital algebra $A$ always admits an enveloping coaction.
\end{prop}

\begin{proof} Let $\prho :A\rightarrow A\otimes H$ be the coaction. As $H$ is
  finite dimensional, its finite dual $H^\circ$ is simply the dual vector
  space $H^*$. The condition that the finite dual separate points is also 
  automatically satisfied in finite dimension. Define the partial action 
of $H^*$ on $A$ by 
\[
f\cdot a=(I\otimes f)\prho (a).
\]
As it came from a partial coaction, it is easy to see that 
$\mbox{dim}(H^* \cdot a)<\infty$, $\forall a\in A$. Choose a basis 
$\{ a_i \}_{i=1}^{n}$ for the subspace $H^* \cdot a$. Then we can prove that 
there are elements $h_i \in H$ for $i=1,\ldots ,n$ such that
\[
f\cdot a =\sum_{i=1}^n \langle h_i , f\rangle a_i .
\]

Let $(B, \varphi)$ be
the standard enveloping action for this partial action on $A$ (remember that
$B\subseteq \mbox{Hom}_k (H^* , A)$). Now, we have to verify
whether the subspace $H^*\triangleright \varphi (a) \subseteq B$ is finite
dimensional. Take $f,g\in H^*$, then
\begin{eqnarray}
(f\triangleright \varphi (a))(g) &=& \varphi (a)(gf) =gf\cdot a=\nonumber\\
&=& \sum_{i=1}^n \langle h_i , gf\rangle a_i =\nonumber\\
&=& \sum_{i=1}^n \langle f\rightharpoonup h_i , g \rangle a_i. \nonumber
\end{eqnarray}
Then, for each $f\in H^*$, we can see that 
\[
f\triangleright \varphi (a)\in \mbox{Hom}_k (H^* ,
\mbox{span}\{ a_i |i=1,\ldots ,n\})\cong
\bigoplus_{i=1}^n H .
\]
As the space $\bigoplus_{i=1}^n H$ is finite dimensional, then
$H\triangleright \varphi (a)$ is also finite dimensional.

Therefore, $B$ is a rational left $H^*$ module, which allows to define a
coaction of $H$ on $B$.
\end{proof}

\section{Partial representations}

Partial representations of groups were first introduced by R. Exel in
\cite{ruy2} and became a powerful tool to investigate the action of semigroups
on algebras. A partial representation of a group $G$ is a map 
$\pi :G\rightarrow A$ on a unital algebra $A$ such that 
\begin{eqnarray}
\label{partialrepgroup}
\mbox{\bf 1)} &\,& \pi (e) = \um , \nonumber\\
\mbox{\bf 2)} &\,& \pi (g)\pi (h) \pi (h^{-1}) = \pi (gh)\pi (h^{-1}), \qquad \forall g,h \in
G, \nonumber\\
\mbox{\bf 3)} &\,& \pi (g^{-1}) \pi (g) \pi (h) = \pi (g^{-1}) \pi (gh),\qquad \forall g,h \in
G.
\end{eqnarray} 
In the reference \cite{dok}, the authors showed that partial actions and
partial representations of a group are intimately related. By one hand, if we
have a partial action of $G$ on a unital algebra $A$ such that each ideal
$D_g$ has unit $1_g$, then there is a partial representation of the group $G$
on the partial crossed product $A\rtimes_{\alpha} G$ given by
\[
\pi (g) =1_g \delta_g .
\]
On the other hand, if there is a partial representation $\pi :G\rightarrow A$,
it is possible to show that $A$ is isomorphic to a partial crossed product
${\bar{A}}\rtimes_{\alpha} G$, where $\bar{A}$ is an abelian subalgebra of
$A$ generated by the elements of the form $\varepsilon_g =\pi (g) \pi (g^{-1})$, 
$\forall g\in G$ and the partial action is given as follows: the ideals $D_g$
are $D_g =\varepsilon_g \bar{A}$ and $\alpha_g :D_{g^{-1}} \rightarrow D_g$
given by $\alpha_g (x)=\pi (g) x \pi (g^{-1})$, $\forall x\in D_{g^{-1}}$.

\begin{prop} Let $\alpha$ be a partial action of $G$ on $A$ such that every idempotent
$1_g$ is central. Then the map $\pi: G \vai End_k(A)$  given by 
\[
\pi (g) (a) = g \cdot a =\alpha_g (a1_{g^{-1}})
\]
defines a partial representation of $G$.
\end{prop}
\begin{proof}

It is easy to see that the first item of (\ref{partialrepgroup})
  holds, because for each $a\in A$
\[
\pi (e)(a)= 1 \cdot a= a. 
\]
Therefore $\pi (g)=I=1_{\mbox{End}_k (A)}$.
As we have shown in the beginning, if the idempotents $1_g = g \cdot \um$ are central then $k \cdot (l \cdot a)= (k \cdot \um)(kl \cdot a) = (kl \cdot a)(k \cdot \um)$. Hence,

\begin{eqnarray}
\pi (g^{-1}) \pi (gh) (a) &=& g^{-1}\cdot (gh \cdot a)=\nonumber\\
& = &  (g^{-1} \cdot \um) (h \cdot a)=\nonumber\\
& = & (g^{-1} \cdot \um) (g^{-1}g \cdot \um) (h \cdot a)=\nonumber\\
& = & (g^{-1} \cdot (g^{-1} \cdot (h \cdot a))) = \nonumber\\
& = &  \pi (g^{-1}) \pi (g) \pi (h) (a). \nonumber
\end{eqnarray}
Since this occurs $\forall a\in A$, we conclude that 
$\pi (g) \pi (h) \pi (h^{-1})=\pi (gh) \pi (h^{-1})$.

In a similar manner, 

\begin{eqnarray}
\pi (gh) \pi (h^{-1}) (a) &=& gh \cdot (h^{-1} \cdot a) =\nonumber\\
& = & (g \cdot a) (gh \cdot \um) =\nonumber\\
& = & (g \cdot \um a ) (gh \cdot \um) =\nonumber\\
& = & (g \cdot \um) (g \cdot  a ) (gh \cdot \um) =\nonumber\\
& = & (g \cdot \um) (ghh^{-1} \cdot  a ) (gh \cdot \um) =\nonumber\\
& = &  (g \cdot \um) (gh \cdot ( h^{-1} \cdot a )) =\nonumber\\
& = &  g \cdot (h \cdot (h^{-1} a))=\nonumber\\
& = &  \pi (g) \pi (h) \pi (h^{-1}) (a)\nonumber
\end{eqnarray}
As this equality holds $\forall a\in A$ then 
$\pi (g^{-1}) \pi (g) \pi (h) (a)=\pi (g^{-1}) \pi (gh) (a)$. Therefore, $\pi$
is a partial representation of $G$.
\end{proof}

Inspired in this previous example, let us try to define a partial
representation of a Hopf algebra $H$. In what follows, we assume that 
$H$ is a Hopf algebra
with invertible antipode (which means that $H$ is Co-Frobenius as a coalgebra, see
for example \cite{romenos}, section 5.4).

\begin{prop} Let $H$ be a Hopf algebra, with invertible antipode, which acts
  partially on a unital algebra $A$. Then the map 
\[
\begin{array}{rccc}
\pi : & H & \rightarrow & \mbox{End}_k (A) , \\
\, & h & \mapsto & \pi (h)
\end{array}
\]
given by $\pi (h) (a)=h\cdot a$, satisfies:
\begin{eqnarray}
\pi (1_H) &=& I, \nonumber\\
\sum \pi (S^{-1} (h_{(2)})) \pi (h_{(1)}) \pi (k) &=& 
\sum \pi (S^{-1} (h_{(2)})) \pi (h_{(1)} k).\nonumber
\end{eqnarray}
\end{prop}

\begin{proof} The first identity is quite obvious. In order to prove the
  second equality, take any alement $a\in A$, then
\begin{eqnarray}
&\, & \sum \pi (S^{-1} (h_{(2)})) \pi (h_{(1)}) \pi (k) (a) = 
\sum S^{-1} (h_{(2)})\cdot (h_{(1)}\cdot (k\cdot a)) =\nonumber\\
&=& \sum (S^{-1} (h_{(3)})\cdot \um)(S^{-1} (h_{(2)})h_{(1)}\cdot (k\cdot a))=
(S^{-1} (h)\cdot \um)(1_H\cdot (k\cdot a))=\nonumber\\
&=& (S^{-1} (h)\cdot \um)(k\cdot a),\nonumber
\end{eqnarray}

On the other hand, we have
\begin{eqnarray}
&\, & \sum \pi (S^{-1} (h_{(2)})) \pi (h_{(1)}k) (a) = 
\sum S^{-1} (h_{(2)})\cdot (h_{(1)}k\cdot a) =\nonumber\\
&=& \sum (S^{-1} (h_{(3)})\cdot \um)(S^{-1} (h_{(2)})h_{(1)}k\cdot a)=
(S^{-1} (h)\cdot \um)(k\cdot a).\nonumber
\end{eqnarray}
\end{proof}

With this result we can propose the following definition:

\begin{defi} Let $H$ be a Hopf algebra with invertible antipode. A partial
  representation of $H$ on a unital algebra $B$ is a linear map
\[
\begin{array}{rccc}
\pi :& H & \rightarrow & B ,\\
\, & h & \mapsto & \pi (h)
\end{array}
\]
such that 
\noindent
\begin{eqnarray}
\label{partialrephopf}
\mbox{\bf 1)} & &  \pi (1_H) =1_B, \\
\mbox{\bf 2)} & &  \sum \pi (S^{-1} (h_{(2)})) \pi (h_{(1)}) \pi (k) = 
\sum \pi (S^{-1} (h_{(2)})) \pi (h_{(1)} k),  \forall h,k \in H. \nonumber
\end{eqnarray}
\end{defi}

Unlike the group case, partial representations of Hopf algebras are not symmetrical with relation to the left and right side. This is not totally unexpected, since we are working with right ideals; we mention, however, that the situation does not seem to improve much if one imposes that $\varphi(A)$ is an ideal in the enveloping action $(B,\varphi)$.

In the group case, it is known that if $G$ acts partially on a unital algebra
$A$ such that every $D_g$ is unital, then there is a partial representation of
$G$ on the partial crossed product. The same occurs with the partial smash
product, that is, if $H$ acts partially on $A$, then there is a partial
representation of $H$ on ${\underline{A\# H}}$.

\begin{teo} Let $H$ be a Hopf algebra with invertible antipode which acts
  partially on a unital algebra $A$. Then the linear map
\[
\begin{array}{rccc}
\pi : & H & \rightarrow & {\underline{A\# H}} , \\
\, & h & \mapsto & \sum (h_{(1)} \cdot \um) \otimes h_{(2)} 
\end{array}
\]
is a partial representation of $H$.
\end{teo}

\begin{proof} The first item in (\ref{partialrephopf}) can be easily derived,
  since
\[
\pi (1_H) =(1_H \cdot \um)\otimes 1_H =\um \otimes 1_H , 
\]
which is the unit in ${\underline{A\# H}}$. 

For the second item in the definition, take $h, k \in H$, then
\begin{eqnarray}
&\, & \sum \pi (S^{-1} (h_{(2)})) \pi (h_{(1)}) \pi (k) =\nonumber\\ 
&=& (\sum (S^{-1} (h_{(4)}) \cdot \um) \otimes S^{-1}(h_{(3)})) 
(\sum (h_{(1)} \cdot \um) \otimes h_{(2)})
(\sum (k_{(1)} \cdot \um) \otimes k_{(2)}) =\nonumber\\
&=& (\sum (S^{-1} (h_{(5)}) \cdot \um)( S^{-1}(h_{(4)})\cdot (h_{(1)} 
\cdot \um)) \otimes S^{-1}(h_{(3)})h_{(2)})
(\sum (k_{(1)} \cdot \um) \otimes k_{(2)}) =\nonumber\\
&=& (\sum (S^{-1} (h_{(3)}) \cdot \um)
( S^{-1}(h_{(2)})\cdot (h_{(1)} \cdot \um)) \otimes 1_H)
(\sum (k_{(1)} \cdot \um) \otimes k_{(2)}) =\nonumber\\
&=& (\sum (S^{-1} (h_{(4)}) \cdot \um)
( S^{-1}(h_{(3)})\cdot \um)( S^{-1}(h_{(2)})h_{(1)} \cdot \um) \otimes 1_H)
(\sum (k_{(1)} \cdot \um) \otimes k_{(2)}) =\nonumber\\
&=& ((S^{-1} (h) \cdot \um) \otimes 1_H)
(\sum (k_{(1)} \cdot \um) \otimes k_{(2)}) =\nonumber\\
&=& \sum (S^{-1} (h) \cdot \um) (k_{(1)} \cdot \um) \otimes k_{(2)}.\nonumber
\end{eqnarray}

On the other hand, we have
\begin{eqnarray}
& \, & \sum \pi (S^{-1} (h_{(2)})) \pi (h_{(1)}k) =\nonumber\\ 
&=& (\sum (S^{-1} (h_{(4)}) \cdot \um) \otimes S^{-1}(h_{(3)})) 
(\sum (h_{(1)}k_{(1)} \cdot \um) \otimes h_{(2)}k_{(2)}) =\nonumber\\
&=& \sum (S^{-1} (h_{(5)}) \cdot \um)( S^{-1}(h_{(4)})\cdot 
 (h_{(1)}k_{(1)} \cdot \um))\otimes  S^{-1}(h_{(3)})h_{(2)}k_{(2)} 
=\nonumber\\
&=& \sum (S^{-1} (h_{(3)}) \cdot \um)( S^{-1}(h_{(2)})\cdot 
(h_{(1)}k_{(1)} \cdot \um))\otimes  k_{(2)} 
=\nonumber\\
&=& \sum (S^{-1} (h_{(4)}) \cdot \um)( (S^{-1}(h_{(3)})\cdot \um) 
(S^{-1}(h_{(2)})h_{(1)}k_{(1)} \cdot \um))\otimes  k_{(2)} 
=\nonumber\\
&=& \sum (S^{-1} (h) \cdot \um) (k_{(1)} \cdot \um) \otimes k_{(2)}.\nonumber
\end{eqnarray}
\end{proof}

\section{Acknowledgements}

The authors would like to thank to Edson R. Álvares, Eduardo Hoefel for fruitfull discussions. The first author (M.M.S.A.) would like to thank Virgínia S. Rodrigues for her fundamental role in establishing the UFPR-UFSC Hopf Seminars. The second author (E.B.) would like to thank the Math Department of UFPR and its staff for their kind hospitality.

\end{document}